\documentclass[11pt]{amsart}
\pdfoutput=1
\usepackage{microtype}
\usepackage{url}
\usepackage[margin=1.45in]{geometry}
\usepackage[all,cmtip,line]{xy}
\usepackage{color}
\usepackage[colorlinks]{hyperref}
\hypersetup{citecolor=blue}

\theoremstyle{definition}
\newtheorem{theorem}{Theorem}

\newtheorem{definition}[theorem]{Definition}
\newtheorem{remark}[theorem]{Remark}

\numberwithin{equation}{section}

\newcommand{\F}{{\mathbb F}}
\newcommand{\R}{{\mathbb R}}
\newcommand{\C}{{\mathbb C}}
\newcommand{\Z}{{\mathbb Z}}
\newcommand{\Q}{{\mathbb Q}}
\newcommand{\Proj}{{\mathbb P}}
\newcommand{\sO}{{\mathcal O}}

\newcommand{\MW}{Mordell-\kern-.16exWeil }
\newcommand{\kbar}{\overline{k}}

\DeclareMathOperator{\NS}{NS}
\DeclareMathOperator{\Pic}{Pic}
\DeclareMathOperator{\Km}{Km}
\DeclareMathOperator{\PGL2}{PGL_2}
\DeclareMathOperator{\I}{I}
\DeclareMathOperator{\II}{II}
\DeclareMathOperator{\III}{III}

\title{$K3$ surfaces associated with curves of genus two}

\author{Abhinav Kumar} 

\thanks{This work is based on the author's doctoral thesis at Harvard
  University, and was supported by a Putnam Fellowship. This paper was
  completed during the author's postdoctoral appointment at Microsoft
  Research.}

\address{Department of Mathematics\\
Massachusetts Institute of Technology\\
Cambridge, MA 02139}
\email{abhinav@math.mit.edu}

\begin{document}

\begin{abstract}
  It is known (\cite{GLD}, \cite{Na}) that there is a unique K3
  surface $X$ which corresponds to a genus $2$ curve $C$ such that $X$
  has a Shioda-Inose structure with quotient birational to the Kummer
  surface of the Jacobian of $C$. In this paper we give an explicit
  realization of $X$ as an elliptic surface over $\Proj^1$ with
  specified singular fibers of type $\II^*$ and $\III^*$. We describe
  how the Weierstrass coefficients are related to the Igusa-Clebsch
  invariants of $C$.
\end{abstract}

\maketitle
\thispagestyle{empty}
\section{Introduction}

In this paper, we study some K3 surfaces of high rank. Precisely, we
consider K3 surfaces which have a Shioda-Inose structure: that is, an
involution which preserves a global $2$-form, such that the quotient
is a Kummer surface. Kummer surfaces are a special class of K3
surfaces which are quotients of abelian surfaces. The Kummer surface
thus carries algebro-geometric information about the abelian
surface. It has N\'{e}ron-Severi rank at least $17$, and so therefore
do the K3 surfaces with Shioda-Inose structure.

These surfaces were first studied by Shioda and Inose \cite{IS}, who
give a description of {\em singular} K3 surfaces, i.e., those with rank
$20$, the maximum possible for a K3 surface over a field of
characteristic zero. They prove that there is a natural one-to-one
correspondence between the set of singular K3 surfaces up to
isomorphism and the set of equivalence classes of positive definite
even integral binary quadratic forms. The result follows that of
Shioda and Mitani \cite{MS} who show that the set of singular abelian
surfaces (that is, those having Picard number $4$) is also in
one-to-one correspondence with the equivalence classes of positive
definite even integral binary quadratic forms. The construction of
Shioda and Inose produces a singular K3 surface by taking a double
cover of a Kummer surface associated to a singular abelian surface and
with a specific type of elliptic fibration. The resulting K3 surface
has an involution such that the quotient is the original Kummer
surface. It also turns out that the lattice of transcendental cycles
on the K3 surface (i.e.\ the orthogonal complement of the
N\'eron-Severi group in the second singular cohomology group) is
isomorphic to the lattice of transcendental cycles on the abelian
surface.

Nikulin studied groups of automorphisms of K3 surfaces in \cite{N1},
and Morrison studied Shioda-Inose structures more extensively in
\cite{M1}, giving necessary and sufficient conditions for a K3 surface
to have a Shioda-Inose structure, in terms of the N\'eron-Severi group
of the K3 surface.

The Shioda-Inose setup takes a K3 surface and produces the data of a
Kummer surface, or of an abelian surface. One may reverse the
question, and ask: given an abelian surface $A$, can we produce a K3
surface with Shioda-Inose structure such that its quotient under the
Nikulin involution is birational to $\Km(A)$?

In particular, when $A$ is the Jacobian of a curve $C$ of genus $2$,
Galluzzi and Lombardo \cite{GLD} show that there is a unique K3
surface (up to isomorphism) with Shioda-Inose structure which
corresponds to $A$ under this construction. Dolgachev, in the appendix
to the paper, proves an isomorphism between the moduli space of
elliptic K3 surfaces with bad fibers of type $E_8$ and $E_7$, and the
moduli space of principally polarized abelian surfaces. However, the
isomorphism is not explicit, and uses a result of Naruki \cite{Na} to
show that the quotient surface constructed is a Kummer surface.

In this paper we describe explicitly the relation between $C$ and the
K3 surface $X$ with $E_8$ and $E_7$ fibers. There exists an
alternative elliptic fibration on $X$, with a $2$-torsion section. The
Nikulin involution on $X$ is translation by the $2$-torsion
section. The quotient elliptic surface under the $2$-isogeny is the
Kummer surface of a unique principally polarized abelian surface
(generically, the Jacobian of a curve of genus $2$). We relate the
Igusa-Clebsch invariants of the genus $2$ curve to the moduli of the
original K3 elliptic surface, thus giving an explicit description of
the map on moduli spaces. Thus, we also partially answer a question of
Kuwata and Shioda \cite{KS}, who ask for explicit elliptic fibrations
on a generic Kummer surface.

This work generalizes a construction of Elkies \cite{E}, who found a
family of K3 elliptic surface with two $E_8$ fibers such that the
quotient is birational to the Kummer surface of a product of two
elliptic curves $E_1$ and $E_2$, and is related to it by a
Galois-invariant map of degree $9$.  Earlier, Inose \cite{In} had
found the isogenous K3 surface to the Kummer surface of $E_1 \times
E_2$, but with maps that were defined over an algebraically closed
ground field. Recent work of Shioda \cite{S2} realizes elliptic K3
surfaces with two $E_8$ fibers as part of a ``Kummer sandwich'' with
the Kummer surface of $E_1 \times E_2$ over an algebraically closed
ground field. Clingher and Doran have also studied the K3 surfaces
corresponding to $E_1 \times E_2$ in \cite{CD}, as well as to $J(C)$
more recently.

Such explicit formulas have some Diophantine applications. Elkies has
used his techniques from \cite{E} to construct elliptic curves of high
\MW rank over $\Q(t)$ and $\Q$. Similarly, he used a K3 surface of
N\'{e}ron-Severi rank $19$ corresponding to a rational point on a
Shimura curve to find an elliptic curve of \MW rank at least $28$ over
$\Q$. An exposition of these techniques may be found in \cite{E2}. As
an application of the results of this paper, we can find explicit
parametrization of some Hilbert modular surfaces. We hope to address
these techniques and results in a future publication.

The outline of the paper is as follows: in Section \ref{sec:Lattices},
we recall facts about integral lattices which appear in this context
as sublattices of N\'{e}ron-Severi groups or cohomology groups of K3
surfaces. In Section \ref{sec:K3}, we cover the necessary background
on K3 surfaces, whereas section \ref{sec:genus2} recalls the theory of
genus $2$ curves. In Section \ref{sec:mainthm}, we state and prove the
main theorem. In section \ref{sec:6points}, we describe a connection
with the moduli space of $6$ points in $\Proj^1$. In section
\ref{sec:product}, we describe how the family of elliptic surfaces
with two $E_8$ fibers may be obtained as a degeneration of the family
of surfaces with $E_8$ and $E_7$ fibers discussed here, and obtain the
compatibility between our formulas and Shioda's.

\section{Preliminaries on lattices} \label{sec:Lattices}

\begin{definition} 
  A {\em lattice} will denote a finitely generated free abelian group
  $\Lambda$ equipped with a symmetric bilinear form $B: \Lambda \times
  \Lambda \rightarrow \R$. An {\em integral lattice} is a lattice
  whose form takes values in $\Z$.
\end{definition}

The {\em signature} of the lattice is the real signature of the form
$B$, written $(r_+, r_-, r_0)$ where $r_+$, $r_-$ and $r_0$ are the
number of positive, negative and zero eigenvalues of $B$, counted with
multiplicity. Usually $r_0$ is omitted if it is zero, i.e., the form
$B$ has zero kernel, in which case we say that $\Lambda$ is {\em
  non-degenerate}. We say $\Lambda$ is {\em even} if $u^2 = B(u,u) \in
2\Z$ for all $u \in \Lambda$. The dual lattice to $\Lambda$ is denoted
$\Lambda^*$. The {\em discriminant} of a non-degenerate lattice is
$|\det(B)| = |\Lambda^*/\Lambda|$. The lattice is said to be {\em
  unimodular} if its discriminant is $1$. Note that in the literature,
the discriminant is frequently defined to be $\det(B) = (-1)^{r_-} |
\Lambda^*/\Lambda|$, in which case a lattice is unimodular iff its
discriminant is $\pm 1$.

For a lattice $\Lambda$ and a real number $\alpha$, we denote by
$\Lambda(\alpha)$ the lattice which has the same underlying group but
with the bilinear form scaled by $\alpha$. The lattice of rank one
with a generator of norm $\alpha$ will be denoted $\langle \alpha
\rangle$.

By a {\em root} of a positive definite lattice, we will mean an
element $x$ such that $x^2 = 2$, whereas for a negative-definite or
indefinite lattice, we will mean an element $x$ such that $x^2 = -2$.

A {\em root lattice} is a lattice that is generated (as an abelian
group) by its roots.  First, let us introduce some familiar root
lattices, through their Dynkin diagrams. The subscript in the name of
the lattice is the dimension of the lattice, which is also the number
of (blank) nodes in the Dynkin diagram. Adding the extra (filled in)
node gives the extended Dynkin diagram, which will arise later in
connection with elliptic surfaces. The labels on the nodes are the
coefficients of a vector in the kernel of the Cartan matrix of the
extended Dynkin diagram.

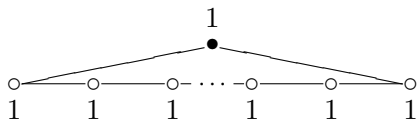
\begin{figure}[h!]
\begin{center}
\begin{picture}(200,50)
\put(10,10){\makebox(0,0){$1$}}
\put(40,10){\makebox(0,0){$1$}}
\put(70,10){\makebox(0,0){$1$}}
\put(100,10){\makebox(0,0){$1$}}
\put(130,10){\makebox(0,0){$1$}}
\put(160,10){\makebox(0,0){$1$}}
\put(85,45){\makebox(0,0){$1$}}

\put(10,20){\circle{4}}
\put(40,20){\circle{4}}
\put(70,20){\circle{4}}
\put(100,20){\circle{4}}
\put(130,20){\circle{4}}
\put(160,20){\circle{4}}
\put(13,20){\line(1,0){24}}
\put(43,20){\line(1,0){24}}
\put(73,20){\line(1,0){4}}
\put(93,20){\line(1,0){4}}
\put(81,20){\circle*{1}}
\put(85,20){\circle*{1}}
\put(89,20){\circle*{1}}
\put(103,20){\line(1,0){24}}
\put(133,20){\line(1,0){24}}
\put(85,35){\circle*{4}}
\put(13,20){\line(5,1){69}}
\put(157,20){\line(-5,1){69}}
\end{picture}
\end{center}
\caption{$A_n  (n\geq 1)$, signature $(n,0)$, discriminant $n+1$.}
\end{figure}

$A_n$ is the positive definite lattice with $n$ generators
$v_1,\ldots, v_n$ with $v_i^2 = 2$ and $v_i \cdot v_j = -1$ if the
vertices $i$ and $j$ are connected by an edge, and $0$ otherwise. It
may be realized as the set of integral points on the hyperplane $\{x
\in \R^{n+1} \big| \sum x_i = 0 \}$.

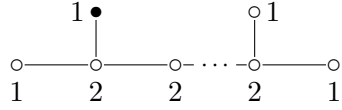
\begin{figure}[h!]
\begin{center}
\begin{picture}(170,40)
\put(10,5){\makebox(0,0){$1$}}
\put(40,5){\makebox(0,0){$2$}}
\put(70,5){\makebox(0,0){$2$}}
\put(100,5){\makebox(0,0){$2$}}
\put(130,5){\makebox(0,0){$1$}}
\put(33,35){\makebox(0,0){$1$}}
\put(107,35){\makebox(0,0){$1$}}

\put(10,15){\circle{4}}
\put(40,15){\circle{4}}
\put(40,35){\circle*{4}}
\put(70,15){\circle{4}}
\put(100,15){\circle{4}}
\put(130,15){\circle{4}}
\put(100,35){\circle{4}}
\put(13,15){\line(1,0){24}}
\put(43,15){\line(1,0){24}}
\put(100,18){\line(0,1){14}}
\put(103,15){\line(1,0){24}}
\put(40,18){\line(0,1){14}}

\put(73,15){\line(1,0){4}}
\put(93,15){\line(1,0){4}}
\put(81,15){\circle*{1}}
\put(85,15){\circle*{1}}
\put(89,15){\circle*{1}}
\end{picture}
\end{center}
\caption{$D_n (n \geq 4)$, signature $(n,0)$, discriminant $4$.}
\end{figure}

$D_n$ can be realized as $\{x \in \Z^n \big| \sum\limits_{i=1}^n x_i
\equiv 0 \bmod 2\}$. It has $2n(n-1)$ roots.

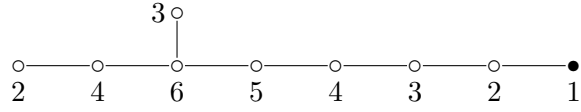
\begin{figure}[h!]
\begin{center}
\begin{picture}(230,40)
\put(10,5){\makebox(0,0){$2$}}
\put(40,5){\makebox(0,0){$4$}}
\put(70,5){\makebox(0,0){$6$}}
\put(100,5){\makebox(0,0){$5$}}
\put(130,5){\makebox(0,0){$4$}}
\put(160,5){\makebox(0,0){$3$}}
\put(190,5){\makebox(0,0){$2$}}
\put(220,5){\makebox(0,0){$1$}}
\put(63,35){\makebox(0,0){$3$}}

\put(10,15){\circle{4}}
\put(40,15){\circle{4}}
\put(70,15){\circle{4}}
\put(100,15){\circle{4}}
\put(130,15){\circle{4}}
\put(70,35){\circle{4}}
\put(160,15){\circle{4}}
\put(190,15){\circle{4}}
\put(220,15){\circle*{4}}

\put(13,15){\line(1,0){24}}
\put(43,15){\line(1,0){24}}
\put(73,15){\line(1,0){24}}
\put(103,15){\line(1,0){24}}
\put(133,15){\line(1,0){24}}
\put(70,18){\line(0,1){14}}
\put(163,15){\line(1,0){24}}
\put(193,15){\line(1,0){24}}

\end{picture}
\end{center}
\caption{$E_8$, signature $(8,0)$, discriminant $1$.}
\end{figure}

One realization of $E_8$ is as the span of $D_8$ and the all-halves
vector $(1/2,\ldots,1/2)$. It has $240$ roots.

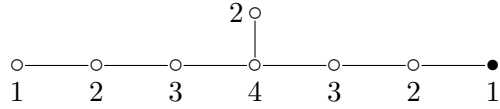
\begin{figure}[h!]
\begin{center}
\begin{picture}(200,40)
\put(10,5){\makebox(0,0){$1$}}
\put(40,5){\makebox(0,0){$2$}}
\put(70,5){\makebox(0,0){$3$}}
\put(100,5){\makebox(0,0){$4$}}
\put(130,5){\makebox(0,0){$3$}}
\put(160,5){\makebox(0,0){$2$}}
\put(190,5){\makebox(0,0){$1$}}
\put(93,35){\makebox(0,0){$2$}}

\put(10,15){\circle{4}}
\put(40,15){\circle{4}}
\put(70,15){\circle{4}}
\put(100,15){\circle{4}}
\put(130,15){\circle{4}}
\put(100,35){\circle{4}}
\put(160,15){\circle{4}}
\put(190,15){\circle*{4}}

\put(13,15){\line(1,0){24}}
\put(43,15){\line(1,0){24}}
\put(73,15){\line(1,0){24}}
\put(103,15){\line(1,0){24}}
\put(133,15){\line(1,0){24}}
\put(163,15){\line(1,0){24}}
\put(100,18){\line(0,1){14}}

\end{picture}
\end{center}
\caption{$E_7$, signature $(7,0)$, discriminant $2$.}
\end{figure}

Taking the orthogonal complement of any root in $E_8$ gives us
$E_7$. It has $126$ roots.

\begin{figure}[h!]
\begin{center}
\begin{picture}(170,60)
\put(10,5){\makebox(0,0){$1$}}
\put(40,5){\makebox(0,0){$2$}}
\put(70,5){\makebox(0,0){$3$}}
\put(100,5){\makebox(0,0){$2$}}
\put(130,5){\makebox(0,0){$1$}}
\put(63,35){\makebox(0,0){$2$}}
\put(63,55){\makebox(0,0){$1$}}

\put(10,15){\circle{4}}
\put(40,15){\circle{4}}
\put(70,15){\circle{4}}
\put(100,15){\circle{4}}
\put(130,15){\circle{4}}
\put(70,35){\circle{4}}
\put(70,55){\circle*{4}}

\put(13,15){\line(1,0){24}}
\put(43,15){\line(1,0){24}}
\put(73,15){\line(1,0){24}}
\put(103,15){\line(1,0){24}}
\put(70,18){\line(0,1){14}}
\put(70,38){\line(0,1){14}}

\end{picture}
\end{center}
\caption{$E_6$, signature $(6,0)$, discriminant $3$.}
\end{figure}
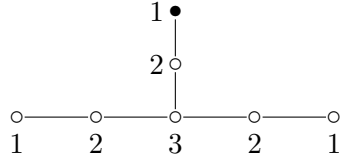

Taking the orthogonal complement of $e_1$ and $e_2$ in $E_8$, where
$e_1,e_2$ are roots such that $e_1 \cdot e_2 = -1$, gives us $E_6$. It
has $72$ roots.

We let the {\em Nikulin lattice} $N$ be the lattice generated by
$v_1,\ldots,v_8$ and $\frac{1}{2}(v_1 + \ldots + v_8)$, with $v_i^2 =
-2$ and $v_i \cdot v_j = 0$ for $i \neq j$. It is isomorphic to
$D_8^*(-2)$.  It has signature $(0,8)$ and discriminant $2^6$. One
checks easily that $N$ has 16 roots, namely $\pm v_i$. In particular,
$N$ is not a root lattice.

Finally, let $U$ be the hyperbolic plane, i.e. the indefinite rank $2$
lattice whose matrix is
\[ \left(\begin{array}{cc}
0 & 1 \\
1 & 0
\end{array}\right). \]

\section{K3 surfaces} \label{sec:K3}

\subsection{Background}

In this section we recall the definition and some basic properties of
K3 surfaces. For more details we refer the reader to \cite{BHPV},
\cite{K} and \cite{M2}. Let $X$ be a smooth projective algebraic
surface over a field $k$.

\begin{definition}
  We say that $X$ is a {\em K3 surface} if $H^1(X,\sO_X) = 0$ and the
  canonical bundle of X is trivial, i.e. $K_X \cong \sO_X$.
\end{definition}

In the sequel, we will assume that $k$ is the field $\C$ of complex
numbers or some subfield of $\C$.

For a K3 surface, one can prove that the middle cohomology $H_X :=
H^2(X,\Z)$ is an even unimodular lattice of signature $(3,19)$,
isomorphic to $E_8(-1)^2\oplus U^3$.

The first Chern class map $H^1(X, \mathcal{O}_X^*) \rightarrow H^2(X,
\Z)$ is injective. Linear equivalence, algebraic equivalence, and
numerical equivalence all coincide for an algebraic K3 surface.

\begin{definition}
  The image of the classes of algebraic divisors $H^1(X,
  \mathcal{O}_X^*)$ in $H^2(X, \Z) \subset H^2(X, \R)$ is a sublattice
  of $H^2(X, \Z)$, which we call the {\em N\'eron-Severi} group of
  $X$, and denote by $\NS(X)$ or $S_X$ or $\Pic(X)$. The orthogonal
  complement of $\NS(X)$ in $H^2(X,\Z)$ is called the {\em
    transcendental} lattice $T_X$ of $X$.
\end{definition}

We have the Hodge decomposition 
$$H^2(X,\Z) \otimes \C  \cong H^{0,2}(X) \oplus H^{1,1}(X) \oplus
H^{2,0}(X) $$ with the vector spaces $H^{0,2}$ and $H^{2,0}$ being
one-dimensional.

\begin{definition} 
  A {\em $(-2)$-curve } on a surface is an irreducible rational curve
  with self-intersection $-2$.
\end{definition}
On a K3 surface, it is an easy exercise using the genus formula to see
that curves of self-intersection $-2$ are exactly the smooth rational
curves.

\subsection{Kummer surfaces}

Here we recall the construction of the Kummer surface associated to an
abelian surface, which is an example of a K3 surface. More details may
be found in \cite{BL},\cite{GH},\cite{H} and \cite{N3}.

Let $A$ be an abelian surface, and let $\iota$ be the involution which
is multiplication by $-1$. Note that $\iota$ fixes any regular
algebraic $2$-form on $A$ (the space of $2$-forms is one-dimensional
since $h^0(X,K) = h^2(X,\sO) = 1$). Now if we blow up the surface $A$
at its sixteen $2$-torsion points, we get a surface $\widetilde{A}$ on
which $\iota$ extends to an involution $\tilde{\iota}$. The quotient
surface $Y = \widetilde{A}/\{1, \tilde{\iota}\}$ is in fact a K3 surface,
and the regular $2$-form is induced from the one on $A$.  The surface
$Y = \Km(A)$ is called the Kummer surface of
$A$. It has sixteen disjoint rational curves $F_1, \ldots, F_{16}$
which are the images of the sixteen exceptional divisors on
$\widetilde{A}$.

The N\'eron-Severi lattice of a Kummer surface has $16$ linearly
independent divisor classes coming from the sixteen rational curves
above. These generate a negative definite lattice, and there is also a
class of a polarization on $\Km(A)$, since it is
projective. Therefore its signature is $(1,r)$ for some $r \geq
16$. In fact, the N\'eron-Severi lattice always contains a particular
lattice of signature $(0,16)$ and discriminant $2^6$, called the
Kummer lattice $K$. We describe its structure.

The set $I := A[2] \cong (\Z/2)^4 \cong \F_2^4$ of 16 elements has a
natural structure of a vector space of dimension $4$ over
$\F_2$. Choose a labeling $I = \{1,2,\ldots,16\}$ (for instance, by
writing $i-1 = \sum_{j=0}^3 b_{ij} 2^j$ for $1 \leq i \leq 16$) and
let $f_1,\ldots,f_{16}$ be the classes of the rational curves
corresponding to the blowups at the $2$-torsion points. Let $Q$ be the
set of $32$ elements consisting of $30$ affine hyperplanes (considered
as subsets of $I$) as well as the empty set and all of $I$. This set
has the structure of a vector space over $\F_2$, the addition
operation being symmetric difference of sets. The set $Q$ is in fact
the Reed-Muller code $\mathcal{R}(1,4)$: the characteristic functions
of the sets in $Q$, viewed as functions from $\F_2^4$ to $\F_2$, are
exactly the polynomials in $4$ variables of degree at most $1$.  For
every $M \in Q$, we have an element $f_M = \frac{1}{2} \sum_{i \in M}
f_i$ of $\sum \Q f_i \in \NS(\Km(A)) \otimes
\Q$. These vectors actually lie in $\NS(\mathop{\rm
  Km}\nolimits(A))$. The lattice generated by the $f_i, (i = 1,\ldots
16)$ and the $f_M$ has discriminant $2^{16}/(2^5)^2 = 2^6$ and it is
called the Kummer lattice.

In section \ref{sec:genus2}, we will say more about the Kummer surface
associated to a principally polarized abelian surface which arises as
the Jacobian of a genus $2$ curve.

\subsection{Shioda-Inose structures}

We now describe the data of a Shioda-Inose structure on a K3 surface.

\begin{definition} 
  An involution $\iota$ on a K3 surface $X$ is called a {\em Nikulin
    involution} if $\iota^*(\omega) = \omega$ for every $\omega \in
  H^{2,0}(X)$ .
\end{definition}

In fact, a Nikulin involution fixes $T_X$ pointwise. Every Nikulin
involution has 8 isolated fixed points. As in the construction of the
Kummer surface, we may blow up these points to get $\widetilde{X}$, which
has eight exceptional curves, and an involution $\tilde{\iota}$. The
quotient $\widetilde{X}/\{1, \tilde{\iota}\}$ of $X$ by a Nikulin
involution is a K3 surface $Y$.

The images of the exceptional curves are $(-2)$-curves on the quotient
K3 surface $Y$. The N\'{e}ron-Severi lattice of $Y$ contains the
Nikulin lattice $N$: it is generated by vectors $c_1,\ldots,c_8$ (the
classes of the eight disjoint $-2$-curves) and $\frac{1}{2} \sum c_i$,
with the form induced by $c_i \cdot c_j = -2\delta_{ij}$.

\begin{definition}
  We say that $X$ admits a {\em Shioda-Inose structure} if there is a
  Nikulin involution $\iota$ on $X$ with rational quotient map $\pi: X
  \rightarrow Y$ such that $Y$ is a Kummer surface and $\pi_*$ induces
  a Hodge isometry $T_X(2) \cong T_Y$.
\end{definition}

If $X$ has a Shioda-Inose structure, let $A$ be the abelian surface
whose Kummer surface is $Y$. Then we have a diagram
$$
\xymatrix{
X \ar@{-->}[dr] & & A \ar@{-->}[dl] \\
& Y &
}
$$
of rational maps of degree $2$, and Hodge isometries $T_X(2) \cong T_Y
\cong T_A(2)$, thus inducing a Hodge isometry $T_X \cong T_A$.

The following theorem of Morrison characterizes K3 surfaces with
Shioda-Inose structures in terms of the N\'{e}ron-Severi lattice.
\begin{theorem}[Morrison  \cite{M1}]
Let $X$ be an algebraic K3 surface. The following are equivalent:
\begin{enumerate}
\item $X$ admits a Shioda-Inose structure.
\item There exists an abelian surface $A$ and a Hodge isometry $T_X \cong T_A$.
\item There is a primitive embedding $T_X \hookrightarrow U^3$.
\item There is an embedding $E_8(-1)^2 \hookrightarrow \NS(X)$.
\end{enumerate}
\end{theorem}

\subsection{Elliptic K3 surfaces}\label{subsec:ellk3}

We recall here a few facts about elliptic surfaces. References for
these are \cite{S1} and \cite{Si}.

\begin{definition}
  An {\em elliptic surface} is a smooth projective algebraic surface $X$
  with a proper morphism $\pi: X \rightarrow C$ to a smooth projective
  algebraic curve $C$, such that
\begin{enumerate}
\item There exists a section $\sigma: C \rightarrow X$.
\item The generic fiber $E$ is an elliptic curve.
\item $\pi$ is relatively minimal.
\end{enumerate}
\end{definition}

Concretely, we will be considering the case $C = \Proj^1$, since a K3
surface cannot have a non-constant map to a curve of positive
genus. We will choose a Weierstrass equation for the generic fiber,
which is an elliptic curve over the function field $\C(\Proj^1) =
\C(t)$, namely
$$
y^2 + a_1(t) xy + a_3(t) y = x^3 + a_2(t) x^2 + a_4(t) x + a_6(t)
$$
where $a_i$ are rational functions of $t$. In fact, by multiplying $x$
and $y$ by suitable rational functions, we can make $a_i(t)$
polynomials in $t$. Furthermore, we may translate $x$ by a rational
function, and $y$ by a $\Q(t)$-linear combination of $1$ and $x$. This
can be done in such a way that the degree of the discriminant is
minimal. We can read out some properties of the surface directly from
the Weierstrass equation. For instance, if $p_a(X)$ is the arithmetic
genus of $X$, then $p_a(X) + 1 = \chi(\mathcal{O}_X)$ is the minimal
$n$ such that $\deg a_i \leq ni$ for $i=1,2,3,4,6$. In particular, for
a K3 surface $X$, we need to have degree $a_i \leq 2i$ (the case $n=1$
corresponds to a rational elliptic surface).

All but finitely many of the fibers of the elliptic surface are
nonsingular and hence elliptic curves. Tate's algorithm \cite{T}
allows us to compute the description of the singular fibers, according
to the Kodaira-N\'{e}ron classification. We note that the reducible
fibers are unions of nonsingular rational curves with
multiplicities. The dual graph of these components is an extended
Dynkin diagram of type $A$, $D$ or $E$. The identity component is
represented by the filled in node in the extended Dynkin diagram,
whereas the multiplicities of the different components of the fiber
are given by the labels on the nodes.

The N\'eron-Severi lattice of $X$ is generated by the classes of all
the sections of $\pi$ (i.e. the \MW group of $X$)
considered as curves on the surface $X$, together with the class $F$
of a fiber, and all the non-identity components of the reducible
fibers. Let $R = \{v \in C(\C)\,|\, F_v \textrm{ is reducible} \}$,
and for each $v \in R$, let $F_v = \pi^{-1}(v) = \Theta_{v,0} +
\sum\limits_{i=1}^{m_v-1} \mu_{v,i} \Theta_{v,i}$, where
$\Theta_{v,0}$ is the component which intersects the identity, and the
other $\Theta_{v,i}$ are the non-identity components. Let $O$ be the
zero section. The intersection pairing satisfies:
\begin{enumerate}
\item for any section $P$, $P^2 = O^2 = - \chi$,
\item $P \cdot F = O \cdot F = 1$,
\item $F^2 = 0$,
\item $O \cdot \Theta_{v,i} = 0$ for $i \geq 1$,
\item $\Theta_{v,i} \cdot \Theta_{w,j} = 0$ for $v \neq w$.
\end{enumerate}
The intersection pairing for $\Theta_{v,i}$ and $\Theta_{v,j}$ is $-2$
if $i = j$, and $0,1,2$ and if $i \neq j$ according to the number of
edges connecting the corresponding nodes in the extended Dynkin
diagram ($2$ occurs only for the types $\I_2$ and
$\III$, which have the extended Dynkin diagram of $A_1$). For
a general section $P$, the intersection pairing with each
$\Theta_{v,i}$ can be computed locally. In particular, for each $v$
exactly one of the intersection numbers is $1$, for some $i$ such that
$\mu_{v,i} = 1$, and the others vanish. The rank of the N\'eron-Severi
group is given by the formula
$$
\rho = r + 2 + \sum_v (m_v -1).
$$
The discriminant of the sublattice $T$ generated by the non-identity
components of all the fibers is $\prod_v m^{(1)}_v$, where $m^{(1)}_v$
is the number of multiplicity one components of $F_v$.

\section{Genus \texorpdfstring{$2$}{2} curves} \label{sec:genus2}
\subsection{Background: Moduli and invariants}

Here we describe the basic geometry and moduli of curves of genus
$2$. For more background we refer the reader to \cite{CF}, \cite{Cl},
\cite{I1}, \cite{Me}. Let $C$ be such a curve defined over a field $k$
of characteristic zero. Then the canonical bundle $K_C$ of $C$ has
degree $2$ and $h^0(C,K_C) = 2$. That is, the corresponding complete
linear system is a $g^1_2$ (and it is the unique $g^1_2$). We
therefore have a map
$$
x: C \rightarrow \Proj^1
$$
which is ramified at $6$ points by the Riemann-Hurwitz formula, and
the function field of $C$ is a quadratic extension of
$k(x)$. Therefore, we may write the equation of $C$ as
$$
y^2 = f(x) = \sum_{i=0}^6 f_i x^i.
$$
The roots of the sextic are the six ramification points of the map $C
\rightarrow \Proj^1$. Their preimages on $C$ are the six Weierstrass
points.

Now, the isomorphism class of $C$ over $\kbar$, the algebraic
closure of $k$, is determined by the isomorphism class of the sextic
$f(x)$, where two sextics are equivalent if there is a transformation
in $\PGL2(\kbar)$ which takes the set of roots (considered inside
$\Proj^1$) to the roots of the other.

Clebsch was the first to determine the invariants of binary
sextics. He defined invariants of $I_2, I_4, I_6, I_{10}$ of weights
$2,4,6,10$ respectively. That is, $I_d$ has degree $d$ in the
coefficients of $f$, and if $f$ transforms to $g$ under the action of
$GL_2(\kbar)$, then there is an element $r \in \kbar$ such that
$I_d(g) = r^d I_d(f)$. Clebsch and Bolza showed that these invariants
determined the sextic up to $\kbar$-equivalence. Therefore, the
point $(I_2(f): I_4(f): I_6(f): I_{10}(f))$ in weighted projective
space determines the isomorphism class of $C$. In fact, $C$ and $C'$
are isomorphic over $k$ iff there is an $r \in k^*$ such that $I_d(f')
= r^d I_d(f)$. Igusa generalized Clebsch's theory to hold in all
characteristics by choosing a different algebraic equation for the
curve $C$ (through an embedding as a quartic in $\Proj^2$ with one
node) and defining invariants $J_2,J_4,J_6,J_8$ and $J_{10}$. He thus
obtained a moduli space of genus two curves defined over $\mathop{\rm
  Spec}\nolimits \Z$. The invariants $I_2,\ldots, I_{10}$ are called
the Igusa-Clebsch invariants whereas $J_2, \ldots, J_{10}$ are called
the Igusa invariants of $C$. We will use the former, since we are
working over a field of characteristic zero.

\begin{remark} 
  If the Igusa-Clebsch invariants of a curve $C$ lie in a field $k$,
  it does not necessarily mean that $C$ can be defined over $k$: there
  is usually an obstruction in $Br_2(k)$. But $C$ can always be
  defined over a quadratic extension of $k$.
\end{remark}

\subsection{Kummer surface}

Let $C$ be a curve of genus $2$, which we can write as
$$
y^2 = f(x) = \sum_{i=0}^6 f_i x^i
$$
Let $\theta_i, i = 1,\ldots, 6$ be the roots of the of the sextic, so that
$$
f(x) = f_6 \prod_{i=1}^6(x-\theta_i).
$$

We shall concern ourselves with the embedding of the singular Kummer
surface as a quartic in $\Proj^3$, which comes from the complete
linear system $2\Theta$, twice the theta divisor which defines the
principal polarization. For a treatment of the quartic surface and the
formulas we use, we refer the reader to \cite{CF},\cite{Co},\cite{H}
and \cite{J}. The quartic is given by the equation
$$
K(z_1,z_2,z_3,z_4) = K_2 z_4^2 + K_1 z_4 + K_0 = 0.
$$
where
\begin{eqnarray*}
K_2 &=& z_2^2-4z_1z_3 ,\cr \cr
K_1 &=& -4z_1^3f_0-2z_1^2z_2f_1-4z_1^2z_3f_2-2z_1z_2z_3f_3-4z_1z_3^2f_4-2z_2z_3^2f_5-4z_3^3f_6, \cr \cr
K_0 &=& -4z_1^4f_0f_2+z_1^4f_1^2-4z_1^3z_2f_0f_3-2z_1^3z_3f_1f_3-4z_1^2z_2^2f_0f_4 \cr
&&+4z_1^2z_2z_3f_0f_5-4z_1^2z_2z_3f_1f_4-4z_1^2z_3^2f_0f_6+2z_1^2z_3^2f_1f_5
\cr
&&-4z_1^2z_3^2f_2f_4+z_1^2z_3^2f_3^2-4z_1z_2^3f_0f_5+8z_1z_2^2z_3f_0f_6
\cr
&& -4z_1z_2^2z_3f_1f_5+4z_1z_2z_3^2f_1f_6-4z_1z_2z_3^2f_2f_5-2z_1z_3^3f_3f_5
\cr
&&-4z_2^4f_0f_6-4z_2^3z_3f_1f_6-4z_2^2z_3^2f_2f_6-4z_2z_3^3f_3f_6-4z_3^4f_4f_6
+z_3^4f_5^2. \cr
\end{eqnarray*}

The $16$ singular points define ordinary double points on the quartic,
which are called nodes. These are given explicitly by the coordinates
\begin{eqnarray*}
p_0 &=& (0:0:0:1) \cr
p_{ij} &=& (1:\theta_i + \theta_j: \theta_i \theta_j: \beta_{0}(i,j)) \cr
\end{eqnarray*}
for $1\leq i < j \leq 6$.

Here $\beta_{0}(i,j)$ is defined as follows.  Let 
$$
f(x) =
(x-\theta_i)(x-\theta_j)h(x),  \textrm{ with } h(x) = \sum_{n=0}^4 h_n
x^n.
$$ Then
$$
\beta_{0}(i,j) = -h_0 - h_2 (\theta_i \theta_j) - h_4 (\theta_i
\theta_j)^2.
$$

The singular point $p_0$ comes from the $0$ point of the Jacobian,
whereas the $p_{ij}$ comes from the $2$-torsion point which is the
difference of divisors $[(\theta_i,0)] - [(\theta_j, 0)]$
corresponding to two distinct Weierstrass points on $C$. The sixteen
singular points are called {\em nodes}.

There are also sixteen hyperplanes in $\Proj^3$ which are tangent to
the Kummer quartic. These are called {\em tropes}. Each trope
intersects the quartic in a conic with multiplicity $2$, and contains
$6$ nodes. Conversely, each node is contained in exactly $6$
tropes. This beautiful configuration is called the $(16,6)$ Kummer
configuration.

The explicit formulas for the tropes are as follows. Six of the tropes
are given by
$$
\theta_i^2 z_1 - \theta_i z_2 + z_3 = 0.
$$
We call this trope $T_i$. It contains the nodes $p_0$ and $p_{ij}$.
The remaining ten tropes are labeled $T_{ijk}$ and correspond to
partitions of $\{1,2,3,4,5,6\}$ into two sets of three, say
$\{i,j,k\}$ and its complement $\{l,m,n\}$. Set
$$
G(X) = (x-\theta_i) (x-\theta_j) (x-\theta_k) = \sum_{r=0}^3 g_r x^r,
$$
$$
H(X) = (x-\theta_l) (x-\theta_m) (x-\theta_n) = \sum_{r=0}^3 h_r x^r.
$$
Then the equation of $T_{ijk}$ is
$$
f_6(g_2 h_0 + g_0 h_2) z_1 + f_6(g_0 + h_0) z_2 + f_6(g_1 + h_1) z_3 + z_4 = 0.
$$

The N\'eron-Severi lattice of the nonsingular Kummer surface contains
classes of rational curves $E_0$ and $E_{ij}$ coming from the nodes,
and $C_i$ and $C_{ijk}$ coming from the tropes. We will denote the
lattice generated by these as $\Lambda_{(16,6)}$. It has signature
$(1,16)$ and discriminant $2^6$ and is the N\'eron-Severi lattice of
the Kummer surface of a generic principally polarized abelian surface.

Let $L$ be the class of a hyperplane section. Furthermore, since
$T_{ijk}$ only depends on the partition $\{1,2,3,4,5,6\} = \{i,j,k\}
\cup \{l,m,n\}$, we may assume $i=1$ and set $C_{jk} = C_{1jk}$ for $1
\leq j < k \leq 6$. Also set $C_{1j} = C_j$ for $1 < j \leq 6$ and
$C_0 = C_1$ to agree with the notation in \cite{Na}. We have the
following intersection numbers and relations in the N\'eron-Severi
lattice.
\begin{eqnarray*}
L^2  &=& 4, \cr
E_0^2 &=& -2, \cr
E_{ij}^2 &=& -2, \cr
E_0 \cdot E_{ij} &=& 0, \cr
E_{ij} \cdot E_{kl} &=& 0  \quad \text{ for } \{i,j\} \neq \{k,l\}, \cr
C_0 &=& \Big(L-E_0-\sum_{k} E_{1k}\Big)/2, \cr
C_{1j} &=& \Big(L-E_0-\sum_{k\neq j} E_{jk}\Big)/2, \cr
C_{jk} &=& (L-E_{ij}-E_{jk}-E_{ik} -E_{lm}-E_{mn}-E_{ln})/2,
\end{eqnarray*}
where $\{l,m,n\}$ is the complementary set to $\{1,j,k\}$.

Projection to a hyperplane from $p_0$ defines a $2$ to $1$ map of the
Kummer to $\Proj^2$, and thus identifies the Kummer surface as a
double cover of $\Proj^2$, ramified along the union of six lines,
which are the projections of the conics $C_i$ (or the tropes
$T_i$). The exchange of sheets gives an involution of the lattice,
which acts by
\begin{eqnarray*}
E_0 &\mapsto& 2L-3E_0, \cr
E_{ij} &\mapsto& E_{ij}, \cr
L &\mapsto& 3L-4E_0.
\end{eqnarray*}

We can explicitly write down the projection to $\Proj^2$ as
$(x_1,x_2,x_3,x_4) \mapsto (x_1,x_2,x_3)$. The involution which is the
exchange of sheets is $(x_1,x_2,x_3,x_4) \mapsto
(x_1,x_2,x_3,-x_4)$. Let $q_0, q_{ij}$ be the projections of the $p_0,
p_{ij}$.

\section{Family of elliptic K3 surfaces associated to genus \texorpdfstring{$2$}{2} curves} \label{sec:mainthm}
\subsection{Elliptic surfaces with \texorpdfstring{$E_8$}{E8} and \texorpdfstring{$E_7$}{E7} fibers}

Let $X$ be an elliptic K3 surface with bad fibers of type $E_8$ and
$E_7$ at $\infty$ and $0$ respectively. A generic such K3 surface has
a N\'eron-Severi lattice $\NS(X) \cong U \oplus E_8(-1) \oplus E_7(-1)$
by Shioda's explicit description of the N\'eron-Severi lattice of an
elliptic surface. This lattice has rank $17$, signature $(1,16)$ and
discriminant $2$. The transcendental lattice $T_X$ has rank $5$,
signature $(2,3)$ and discriminant $2$. We deduce that $T_X \cong U^2
\oplus \langle -2 \rangle$.

The transcendental lattice of a generic principally polarized abelian
surface, that is, the Jacobian $J(C)$ for $C$ a generic curve of genus
$2$, also satisfies the same property, since the N\'eron-Severi of
$J(C)$ is generated by the theta divisor, which has self-intersection
$2$ by the genus formula on the abelian surface
$$
2 = 2g -2 = C\cdot(C+K)= C^2.
$$
Therefore the orthogonal complement in $H^2(J(C),\Z) \cong U^3$ is
exactly $U^2 \oplus \langle -2\rangle $. We expect that the elliptic
K3 surface $X$ has a Shioda-Inose structure such that the quotient by
the Nikulin involution gives the Kummer surface of a principally
polarized abelian surface $\Km(J(C))$. Galluzzi,
Lombardo and Dolgachev \cite{GLD} prove that, in fact, $X$ corresponds
to a unique $C$ up to isomorphism. However, an explicit identification
of the quotient as a Kummer surface was not known. Below, we give an
explicit construction of the correspondence.

We begin with the K3 surface $X$ given by the equation
$$
y^2 = x^3 + t^3(at + a')x + t^5(b''t^2 + bt + b').
$$
It is easily checked that the surface $X$ has an $\II^*$ or
$E_8$ fiber at $t = \infty$ and a $\III^*$ or $E_7$ fiber at
$t = 0$. Now we describe another elliptic fibration on $X$.

Applying the transformation $(x,y,t) = (x't'^2/b''^2,y't'^2/b''^3,
t'/b'')$, we get
$$
y'^2 = t'^3 + (x'^3 + ax' + b)t'^2 + b''(a'x'+b') t'
$$
and again replacing $(x',y',t')$ by $(x,y,t)$ for convenience gives
finally
$$
y^2 = t^3 + (x^3 + ax + b)t^2 + b''(a'x+b') t
$$
which is an elliptic surface over the $x$-line with an
$\I^*_{10}$ or $D_{14}$ fiber at $x = \infty$, an
$\I_2$ or $A_1$ fiber at $x = -b'/a'$ and a $2$-torsion
section $(y,t) = (0,0)$.  The translation by the $2$-torsion section
is a Nikulin involution. We write down the isogenous elliptic surface
$Y$ as
$$
y^2 = t^3 -2(x^3 + ax + b)t^2 + ((x^3 + ax + b)^2 -  4b''(a'x+b'))t.
$$
This is an elliptic surface over the $x$-line with an $\I^*_5$
or $D_9$ fiber at $x = \infty$ and $\I_2$ or $A_1$ fibers at
the roots of the sextic $(x^3 + ax + b)^2 - 4 b''(a'x+b')$, and with a
$2$-torsion section $(t,y) = (0,0)$. The N\'eron-Severi lattice of a
generic such surface has signature $(1,16)$ and discriminant $4\cdot
2^6/2^2 = 2^6$. In fact, we will identify it with the N\'eron-Severi
lattice of a generic Kummer surface (which we call the $(16,6)$
lattice) in a later section. This will lead to the identification of
the Kummer surface of $J(C)$ as an elliptic K3 surface with a bad
fiber of type $\I^*_5$ at $\infty$, with $\I_2$ fibers
at the roots of a sextic derived from $C$, and with a $2$-torsion
section.

\begin{remark}
  If $E$ is an elliptic curve with a $2$-torsion point, written in the
  Weierstrass form
$$
y^2 = x^3 + ax^2 + bx
$$
with a $2$-torsion point $P = (0,0)$, then the $2$-isogenous curve $E'
= E/\{O,P\}$ is given by
$$
Y^2 = X^3 - 2aX^2 + (a^2-4b)X.
$$
The isogeny $\phi: E \rightarrow E'$ is given by
$$
(x,y) \rightarrow \bigg(\frac{y^2}{x^2}, \frac{y(b-x^2)}{x^2}\bigg).
$$
The dual isogeny $\hat{\phi}$ is given by
$$
(X,Y) \rightarrow \bigg(\frac{Y^2}{4X^2}, \frac{Y(a^2-4b-X^2)}{8X^2}\bigg).
$$
\end{remark}

\subsection{Main theorem}

In this section, we state the main theorem, which gives the parameters
$(a,a',b,b',b'')$ of the K3 elliptic surface with $E_8$ and $E_7$
fibers (so far conjecturally) associated to a genus $2$ curve $C$ in
terms of the Igusa-Clebsch invariants of $C$. The proof will be given
in the following sections.

\begin{theorem}
  Let $C$ be a curve of genus two, and $Y = \mathop{\rm
    Km}\nolimits(J(C))$ the Kummer surface of its Jacobian. Let $I_2,
  I_4, I_6, I_{10}$ be the Igusa-Clebsch invariants of $Y$. Then there
  is an elliptic fibration on $Y$ for which the Weierstrass equation
  may be written
$$
y^2 = x^3 -2\bigg(t^3 - \frac{I_4}{12} t + \frac{I_2 I_4 - 3 I_6}{108}\bigg)x^2 + \bigg(\bigg(t^3 - \frac{I_4}{12} t + \frac{I_2 I_4 - 3 I_6}{108}\bigg)^2 + I_{10}\bigg(t - \frac{I_2}{24}\bigg)\bigg) x.
$$
There is an elliptic K3 surface $X$ given by
$$
y^2 = x^3 - t^3\bigg( \frac{I_4}{12} t +1\bigg)x + t^5\bigg(\frac{I_{10}}{4} t^2 + \frac{I_2 I_4 - 3 I_6}{108} t + \frac{I_2}{24}\bigg)
$$
with fibers of type $E_8$ and $E_7$ at $t = \infty$ and $t=0$
respectively, and a Nikulin involution on $X$, such that the quotient
K3 surface is $Y$.
\end{theorem}

\begin{remark}
  Note that the correspondence of the K3 surface $X$ with the genus
  $2$ curve $C$ is Galois invariant, i.e. $X$ is defined over the
  field of definition of $C$. The proof of the theorem will involve
  making some non-Galois invariant choices (for the level $2$
  structure of the Jacobian of $C$) but we will show that $X$ as well
  as the Shioda-Inose structure on $X$ is independent of these
  choices.
\end{remark}

\begin{remark}
The Nikulin involution on $X$ may be written as follows:
$$
(x,y,t) \mapsto \bigg(\frac{16\, x\,(-x+I_2 t^2/24)^2}{I_{10}^2 \,\,t^8} \, , \,\frac{-64\, y\,(-x+I_2 t^2/24)^3}{I_{10}^3 \,\,t^{12}} \, ,\, \frac{4\,(-x+I_2 t^2/24)}{I_{10} \,\,t^3}   \bigg).
$$

\end{remark}

\subsection{N\'eron-Severi lattices}

In this and the next section, we give the details of how to put an
elliptic fibration on the Kummer surface of a Jacobian of a curve of
genus $2$, with a $2$-torsion section, a $\I^*_5$ fiber and
six $\I_2$ fibers. We first identify the N\'eron-Severi
lattices involved (namely $\Lambda_{(16,6)}$ and $(D_9 \oplus A_1^6
\oplus U)^+$).  We use the work of Naruki \cite{Na}, which gives an
embedding of the lattice $N \oplus E_8(-1)$ inside $\Lambda_{(16,6)}$.

First, we start with the N\'eron-Severi lattice of the K3 surface $X$
which has $E_8$ and $E_7$ fibers. The roots of the $\NS(X)$ which
correspond to the smooth rational curves on $X$ are drawn below (we
use the notation from \cite{GLD}). These are in fact all the roots
corresponding to smooth rational curves on $X$, by a result of Nikulin
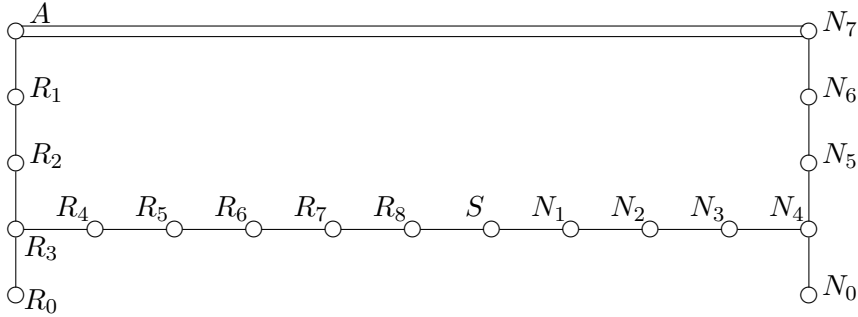
\cite{N2}.

\begin{figure}[h!]
\begin{center}
\begin{picture}(350,150)
\put(10,10){\circle{6}}
\put(10,35){\circle{6}}
\put(10,60){\circle{6}}
\put(10,85){\circle{6}}
\put(10,110){\circle{6}}
\put(310,10){\circle{6}}
\put(310,35){\circle{6}}
\put(310,60){\circle{6}}
\put(310,85){\circle{6}}
\put(310,110){\circle{6}}

\put(40,35){\circle{6}}
\put(70,35){\circle{6}}
\put(100,35){\circle{6}}
\put(130,35){\circle{6}}
\put(160,35){\circle{6}}
\put(190,35){\circle{6}}
\put(220,35){\circle{6}}
\put(250,35){\circle{6}}
\put(280,35){\circle{6}}

\put(12,112){\line(1,0){296}}
\put(12,108){\line(1,0){296}}
\put(10,13){\line(0,1){19}}
\put(10,38){\line(0,1){19}}
\put(10,63){\line(0,1){19}}
\put(10,88){\line(0,1){19}}

\put(310,13){\line(0,1){19}}
\put(310,38){\line(0,1){19}}
\put(310,63){\line(0,1){19}}
\put(310,88){\line(0,1){19}}

\put(13,35){\line(1,0){24}}
\put(43,35){\line(1,0){24}}
\put(73,35){\line(1,0){24}}
\put(103,35){\line(1,0){24}}
\put(133,35){\line(1,0){24}}
\put(163,35){\line(1,0){24}}
\put(193,35){\line(1,0){24}}
\put(223,35){\line(1,0){24}}
\put(253,35){\line(1,0){24}}
\put(283,35){\line(1,0){24}}

\put(13,5){$R_0$}
\put(13,25){$R_3$}
\put(15,60){$R_2$}
\put(15,85){$R_1$}
\put(15,113){$A$}
\put(25,40){$R_4$}
\put(55,40){$R_5$}
\put(85,40){$R_6$}
\put(115,40){$R_7$}
\put(145,40){$R_8$}
\put(180,40){$S$}
\put(205,40){$N_1$}
\put(235,40){$N_2$}
\put(265,40){$N_3$}
\put(295,40){$N_4$}
\put(315,10){$N_0$}
\put(315,60){$N_5$}
\put(315,85){$N_6$}
\put(315,110){$N_7$}

\end{picture}
\end{center}
\end{figure}

There is an elliptic fibration on $X$ which has $R_8 + 2R_7 + 3R_6 +
4R_5 + 5R_4 + 6 R_3 + 4 R_2 + 2 R_1 + 3 R_0$ as an $\II^*$ or
$E_8$ fiber, $N_7 + 2N_6 + 3 N_5 + 4 N_4 + 3N_3 + 2 N_2 + N_1 + 2N_0$
as a $\III^*$ or $E_7$ fiber, and $S$ as the zero
section. This is the fibration over $\Proj^1_t$. The fibration over
$\Proj^1_x$ has the $\I^*_{10}$ or $D_{14}$ fiber given by
$R_0 + R_2 + 2( R_3 + R_4 + R_5 + R_6 + R_7 + R_8 + S + N_1 + N_2 +
N_3 + N_4) + N_0 + N_5$, an $\I_2$ or $A_1$ fiber $A + N_7$, a
$2$-torsion section (say $R_1$) and a zero section $N_6$.

The Nikulin involution $\sigma$ is translation by the $2$-torsion
section. It reflects the above picture about its vertical axis of
symmetry. There are two obvious copies of $E_8(-1)$ switched by
$\sigma$, namely the sublattices of $\NS(X)$ generated by the roots
$\{S, N_1, N_2, N_3, N_4, N_0, N_5, N_6\}$ and $\{R_7, R_6, R_5, R_4,
R_3, R_0, R_2, R_1\}$. Next, we write down some roots on $\NS(Y)$,
where $Y$ is the quotient K3 surface of $X$ by the involution. As we
have described, $Y$ has six $\I_2$ or $A_1$ fibers
$Q_{13}+Q_{14},\ldots, Q_{23}+Q_{24}$, a $\I^*_5$ or $D_9$
fiber, namely $Q_1 + Q_2 + 2(Q_3+Q_4+Q_5 + Q_6+Q_7+ Q_8) + Q_9 +
Q_{10}$, with zero section $O = Q_{11}$ and a $2$-torsion section $T =
Q_{12}$.

\begin{figure}[h!]
\begin{center}
\begin{picture}(350,170)
\put(10,10){\circle{6}}
\put(10,160){\circle{6}}
\put(160,10){\circle{6}}
\put(160,160){\circle{6}}
\put(10,80){\circle{6}}
\put(40,80){\circle{6}}
\put(70,80){\circle{6}}
\put(100,80){\circle{6}}
\put(130,80){\circle{6}}
\put(160,80){\circle{6}}

\put(260,10){\circle{6}}
\put(260,160){\circle{6}}

\put(200,100){\circle{6}}
\put(200,70){\circle{6}}

\put(230,100){\circle{6}}
\put(230,70){\circle{6}}

\put(260,100){\circle{6}}
\put(260,70){\circle{6}}

\put(290,100){\circle{6}}
\put(290,70){\circle{6}}

\put(320,100){\circle{6}}
\put(320,70){\circle{6}}

\put(350,100){\circle{6}}
\put(350,70){\circle{6}}

\put(202,102){\line(1,1){56}}
\put(202,68){\line(1,-1){56}}

\put(231,102){\line(1,2){28}}
\put(231,68){\line(1,-2){28}}

\put(260,103){\line(0,1){54}}
\put(260,67){\line(0,-1){54}}

\put(289,102){\line(-1,2){28}}
\put(289,68){\line(-1,-2){28}}

\put(318,102){\line(-1,1){56}}
\put(318,68){\line(-1,-1){56}}

\put(348,102){\line(-3,2){85}}
\put(348,68){\line(-3,-2){85}}

\put(163,10){\line(1,0){94}}
\put(163,160){\line(1,0){94}}

\put(13,80){\line(1,0){24}}
\put(43,80){\line(1,0){24}}
\put(73,80){\line(1,0){24}}
\put(103,80){\line(1,0){24}}
\put(133,80){\line(1,0){24}}

\put(10,13){\line(0,1){64}}
\put(160,13){\line(0,1){64}}
\put(10,83){\line(0,1){74}}
\put(160,83){\line(0,1){74}}

\put(198,72){\line(0,1){26}}
\put(202,72){\line(0,1){26}}

\put(228,72){\line(0,1){26}}
\put(232,72){\line(0,1){26}}

\put(258,72){\line(0,1){26}}
\put(262,72){\line(0,1){26}}

\put(288,72){\line(0,1){26}}
\put(292,72){\line(0,1){26}}

\put(318,72){\line(0,1){26}}
\put(322,72){\line(0,1){26}}

\put(348,72){\line(0,1){26}}
\put(352,72){\line(0,1){26}}

\put(15,10){$Q_2$}
\put(15,160){$Q_1$}
\put(15,87){$Q_3$}
\put(35,87){$Q_4$}
\put(65,87){$Q_5$}
\put(95,87){$Q_6$}
\put(125,87){$Q_7$}
\put(145,87){$Q_8$}
\put(145,10){$Q_9$}
\put(140,160){$Q_{10}$}
\put(265,5){$Q_{11} = O$}
\put(265,160){$Q_{12} =  T$}

\put(180,70){$Q_{13}$}
\put(210,70){$Q_{15}$}
\put(240,70){$Q_{17}$}
\put(270,70){$Q_{19}$}
\put(300,70){$Q_{21}$}
\put(330,70){$Q_{23}$}

\put(180,100){$Q_{14}$}
\put(210,100){$Q_{16}$}
\put(240,100){$Q_{18}$}
\put(270,100){$Q_{20}$}
\put(300,100){$Q_{22}$}
\put(330,100){$Q_{24}$}

\end{picture}
\end{center}
\end{figure}

It is easily checked that the rational components of the $E_8$ fiber
described above map as following: $N_6 \mapsto O$ (recall that $N_6$
is the zero section of the $D_{14}$ fibration, on which the quotient
map is an isogeny of elliptic surfaces), $N_5 \mapsto Q_9, N_4 \mapsto
Q_8, N_0 \mapsto Q_{10}, N_3 \mapsto Q_7, N_2 \mapsto Q_6, N_1 \mapsto
Q_5, S \mapsto Q_4$. Hence, we see a natural copy of $E_8$ within the
N\'eron-Severi lattice of $Y$. On the other hand, we can also see
eight roots orthogonal to all the generators of $E_8$ as well as to
each other, namely $Q_{14}, Q_{16}, Q_{18}, Q_{20}, Q_{22}, Q_{24},
Q_1$ and $Q_2$.

Now, we use the calculation of \cite{Na} which gives an explicit
embedding of $N \oplus E_8(-1)$ inside the N\'eron-Severi lattice of a
Kummer surface of a generic principally polarized abelian surface, or
$\Lambda_{(16,6)}$. We extend this embedding to get an identification
of $\Lambda_{(16,6)}$ with $\NS(Y)$, i.e. the lattice generated by the
roots in the diagram above.

The identification is as follows:

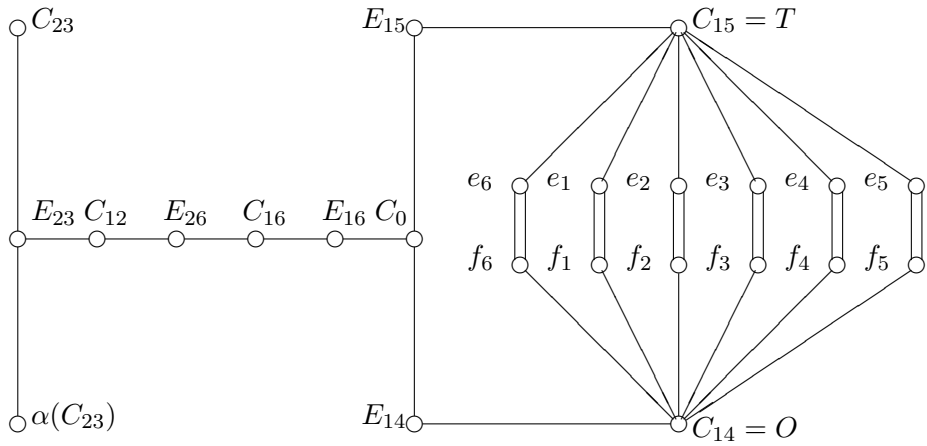
\begin{figure}[h!]
\begin{center}
\begin{picture}(350,170)
\put(10,10){\circle{6}}
\put(10,160){\circle{6}}
\put(160,10){\circle{6}}
\put(160,160){\circle{6}}
\put(10,80){\circle{6}}
\put(40,80){\circle{6}}
\put(70,80){\circle{6}}
\put(100,80){\circle{6}}
\put(130,80){\circle{6}}
\put(160,80){\circle{6}}

\put(260,10){\circle{6}}
\put(260,160){\circle{6}}

\put(200,100){\circle{6}}
\put(200,70){\circle{6}}

\put(230,100){\circle{6}}
\put(230,70){\circle{6}}

\put(260,100){\circle{6}}
\put(260,70){\circle{6}}

\put(290,100){\circle{6}}
\put(290,70){\circle{6}}

\put(320,100){\circle{6}}
\put(320,70){\circle{6}}

\put(350,100){\circle{6}}
\put(350,70){\circle{6}}

\put(202,102){\line(1,1){56}}
\put(202,68){\line(1,-1){56}}

\put(231,102){\line(1,2){28}}
\put(231,68){\line(1,-2){28}}

\put(260,103){\line(0,1){54}}
\put(260,67){\line(0,-1){54}}

\put(289,102){\line(-1,2){28}}
\put(289,68){\line(-1,-2){28}}

\put(318,102){\line(-1,1){56}}
\put(318,68){\line(-1,-1){56}}

\put(348,102){\line(-3,2){85}}
\put(348,68){\line(-3,-2){85}}

\put(163,10){\line(1,0){94}}
\put(163,160){\line(1,0){94}}

\put(13,80){\line(1,0){24}}
\put(43,80){\line(1,0){24}}
\put(73,80){\line(1,0){24}}
\put(103,80){\line(1,0){24}}
\put(133,80){\line(1,0){24}}

\put(10,13){\line(0,1){64}}
\put(160,13){\line(0,1){64}}
\put(10,83){\line(0,1){74}}
\put(160,83){\line(0,1){74}}

\put(198,72){\line(0,1){26}}
\put(202,72){\line(0,1){26}}

\put(228,72){\line(0,1){26}}
\put(232,72){\line(0,1){26}}

\put(258,72){\line(0,1){26}}
\put(262,72){\line(0,1){26}}

\put(288,72){\line(0,1){26}}
\put(292,72){\line(0,1){26}}

\put(318,72){\line(0,1){26}}
\put(322,72){\line(0,1){26}}

\put(348,72){\line(0,1){26}}
\put(352,72){\line(0,1){26}}

\put(15,10){$\alpha(C_{23})$}
\put(15,160){$C_{23}$}
\put(15,87){$E_{23}$}
\put(35,87){$C_{12}$}
\put(65,87){$E_{26}$}
\put(95,87){$C_{16}$}
\put(125,87){$E_{16}$}
\put(145,87){$C_0$}
\put(140,10){$E_{14}$}
\put(140,160){$E_{15}$}
\put(265,5){$C_{14} = O$}
\put(265,160){$C_{15} =  T$}

\put(180,70){$f_6$}
\put(210,70){$f_1$}
\put(240,70){$f_2$}
\put(270,70){$f_3$}
\put(300,70){$f_4$}
\put(330,70){$f_5$}

\put(180,100){$e_6$}
\put(210,100){$e_1$}
\put(240,100){$e_2$}
\put(270,100){$e_3$}
\put(300,100){$e_4$}
\put(330,100){$e_5$}
\end{picture}
\end{center}
\end{figure}

Here $\alpha(C_{23}) = C_{23} + L - 2E_0$.

The class of the fiber is
\begin{eqnarray*}
F &=& C_{23} + \alpha(C_{23}) + 2(E_{23} + C_{12} + E_{26} + C_{16} + E_{16} + C_0) + E_{15} + E_{14}  \cr
&=& 5(L-E_0)-3E_{12} - 2(E_{13} + E_{46} + E_{56}) - (E_{24} + E_{25} + E_{36} + E_{45})
\end{eqnarray*}
and $e_1, \ldots, e_6, f_1, \ldots, f_6$ are given by
\begin{eqnarray*}
e_1 &=& (L-E_0) - (E_{12} + E_{46}) \cr
e_2 &=& 2(L-E_0) - (E_{12} + E_{13} + E_{24} + E_{46} + E_{56}) \cr
e_3 &=& 3(L-E_0) - 2E_{12} - (E_{13} + E_{24} + E_{36} + E_{45} + E_{46} + E_{56}) \cr
e_4 &=& 4(L-E_0) - 2(E_{12} + E_{13} + E_{46}) - (E_{24} + E_{25} + E_{36} + E_{45} + E_{56}) \cr
e_5 &=& 5(L-E_0) - 3E_{12} -2(E_{13} + E_{46} + E_{56}) - (E_{24} + E_{25} + E_{34} + E_{36} + E_{45}) \cr
e_6 &=& E_{35} \cr
(f_i &=& F-e_i \textrm{  for all }i) \cr
f_4 &=& (L-E_0) - (E_{12} + E_{56}) \cr
f_3 &=& 2(L-E_0) - (E_{12} + E_{13} + E_{25} + E_{46} + E_{56}) \cr
f_2 &=& 3(L-E_0) - 2E_{12} - (E_{13} + E_{25} + E_{36} + E_{45} + E_{46} + E_{56}) \cr
f_1 &=& 4(L-E_0) - 2(E_{12} + E_{13} + E_{56}) - (E_{24} + E_{25} + E_{36} + E_{45} + E_{46}) \cr
f_6 &=& 5(L-E_0) - 3E_{12} -2(E_{13} + E_{46} + E_{56}) - (E_{24} + E_{25} + E_{34} + E_{36} + E_{45}) \cr
f_5 &=& E_{34}. \cr
\end{eqnarray*}

Notice that under the simple transposition $(4 5)$ of indices we have
the permutation of fibers $\tau = (14)(23)(56)$ and in fact $e_i
\mapsto f_{\tau(i)}$, $f_i \mapsto e_{\tau(i)}$.

\subsection{Completion of proof}

Next, we describe how to use all this information from the
N\'eron-Severi group to construct $x$, $y$ and $t$ in the Weierstrass
equation for $Y = \Km(J(C))$:
$$
y^2 = x^3 + a(t) x^2 + b(t) x.
$$
Consider the class of the fiber $F \in \NS(\mathop{\rm
  Km}\nolimits(J(C)))$, given by
$$
F = 5(L-E_0)-3E_{12} - 2(E_{13} + E_{46} + E_{56}) - (E_{24} + E_{25} + E_{36} + E_{45}).
$$
We can write down the parameter on the base by computing explicitly
the sections of $H^0(Y,\mathcal{O}_Y(F))$. This linear system consists
of (the pullback of) quintics passing through the points $q_0$ and
$q_{ij}$ which pass through $q_{24},q_{25},q_{36}, q_{45}$, having a
double point at $q_{13}, q_{46}, q_{56}$ and a triple point at
$q_{12}$. This linear system is two-dimensional, and taking the ratio
of two linearly independent sections gives us the parameter $t$ on the
base, $\Proj^1$, for the elliptic fibration.  Now, $t$ is only
determined up to the action of $\PGL2$, but the first restriction we
make is to put the $\I^*_5$ fiber at $t = \infty$, which fixes
$t$ up to affine linear transformations. Any elliptic K3 surface with
a $2$-torsion section can be written in the form
$$
y^2 = x^3 -2 q(t) x^2 + p(t) x
$$
with $p(t)$ of degree at most $8$ and $q(t)$ of degree at most
$4$. The $2$-torsion section is $(x,y) = (0,0)$. The discriminant of
this elliptic surface is a multiple of $p^2(q^2 - p)$.  In fact, we
see that $p$ must have degree exactly $6$, because there are six
$\I_2$ fibers that the zero and $2$-torsion sections meet in
different components. The positions $t_1,\ldots,t_6$ of the
$\I_2$ fibers are the roots of the polynomial $p(t) = p_0
\prod_{i=1}^6 (t-t_i)$.  Now $t$ is determined up to transformations
of the form $t \mapsto at+b$. To have exactly a $\I_5^*$ fiber
at $\infty$, we must have $p(t) = q(t)^2 + r(t)$ where $q(t)$ is a
monic cubic polynomial and $r(t)$ is a {\em linear} polynomial in
$t$. We can further fix $t$ up to scalings $t \mapsto at$ by
translating $t$ so that the quadratic term of $q(t)$ vanishes. We
notice that the top coefficient $p_0$ of $p(t)$ is a square, and so by
scaling $t, x, y$ appropriately, we may assume $p_0 = 1$, i.e. that
$p(t)$ and $q(t)$ are monic.

Now we describe how to obtain $x$. It is a Weil function, so that the
horizontal component of its divisor equals $2T - 2O$, and the vertical
component is uniquely determined by that fact that $(x)$ is linearly
(and hence numerically) equivalent to zero. So we deduce that the
divisor of $x$ is $2T - 2O + Q_{10} -Q_9 + Q_{14} + Q_{16} + Q_{18} +
Q_{20} + Q_{22} + Q_{24} - 3F_0$, where
$$F_0 = Q_1 + Q_2 + 2(Q_3+Q_4+Q_5 + Q_6+Q_7+ Q_8) + Q_9 + Q_{10}$$ is the $D_9$ fiber.

To convert this to formulas, we compute the functions which cut out
$Q_{16},\ldots,Q_{24}$, $F_0$, $T$ and $O$. There is a quintic $s_1$
which cuts out $O = C_{14}$. Now, notice that the $D_9$ fiber contains
$C_{12}, C_{16}$ and $C_0$. Therefore $s_1$ is divisible by $T_2, T_6$
and $T_1$. We write
$$
s_1 = q_1 T_1 T_2 T_6
$$
with a quadratic $q_1$. Next, we know that $T_4$ cuts out $C_{14} = O$
and $T_5$ cuts out $C_{15} = T$. To find, for instance, the function
which cuts out $e_2$, we find the quadratic (unique up to constants)
which passes through $q_{12}, q_{13}, q_{24}, q_{46}, q_{56}$. Call
this function $e_2$, by abuse of notation. Similarly, we find
$e_1,\ldots, e_5$. We also note that the factor of $T_5$ in the
numerator of $x$, which gives a zero along $T$, also gives a zero
along $e_6 = E_{35}$ owing to the fact that $T = C_{15}$ intersects
$E_{35}$ (recall that we are working with the singular Kummer surface,
on which the image of the curve $E_{35}$ is just a single point).

Putting everything together, we can write $x$ up to scaling as a
quotient of two homogeneous polynomials of degree $16$ as follows:
$$
x = \frac{e_1 e_2 e_3 e_4 e_5 T_5}{s_1^3 T_4} = \frac{e_1 e_2 e_3 e_4 e_5 T_5}{ (T_1 T_2 T_6 q_1)^3 T_4}
$$
Finally, we have to scale $x$ and $t$ so that $x^3 + a(t) x^2 + b(t)
x$ becomes a square of a function $y$ on the Kummer.

We note that in the equation of the Kummer
$$
K_2z_4^2 + K_1 z_4 + K_0 = 0
$$
we can complete the square for $z_4$ to obtain
$$
(K_2 z_4 + K_1/2)^2 = K_1^2/4 - K_0 K_2 = 4 T_1 T_2 T_3 T_4 T_5 T_6
$$
We let $y$ be a constant multiple of
$$
\frac{e_1 e_2 e_3 e_4 e_5 (K_2 z_4 + K_1/2)}{T_1^5 T_2^3 T_4^2 T_6^4 q_1^2},
$$
a quotient of two homogeneous polynomials of degree $18$, and verify
that this makes the Weierstrass equation hold. The computation is
carried out in a Maxima program, which is available from the arXiv.org
e-print archive. This paper is available as
\href{http://arxiv.org/abs/math/0701669}{\texttt{math.AG/0701669}}. To
access the auxiliary file, download the source file for the
paper. That will produce not only the \LaTeX \, files for the paper
but also the computer algebra code. The code takes about half an hour
to run on a $2.1$ GHz computer.

We noted earlier that the permutation $(45)$ on the $A_1$ fibers by
$\tau = (14)(23)(56)$, and takes $e_i$ to $f_{\tau(i)}$. That is, it
switches the components intersecting the identity and $2$-torsion
sections as well. In addition, it switches the zero section $C_{14}$
and the $2$-torsion section $C_{15}$, and on the $D_9$ fibers it
switches the two near leaves $E_{15}$ and $E_{14}$, namely, again the
components intersecting $T$ and $O$. On the other hand, consider the
action on $\NS(Y)$ induced by the translation by $T$. Under this map,
$T$ and $O$ get swapped, the $2$-torsion and identity components of
the $D_9$ and $A_1$ fibers all get switched, and the far leaves of the
$D_9$ fiber also get switched (this can be seen, for instance, from
the fact that the simple components of the special $D_9$ fiber form a
group compatible with the group law on the generic elliptic
curve). The locations of the $A_1$ fibers themselves are
fixed. Therefore the effect of the permutation $45$ is the same as
translation by $2$-torsion composed with a pure involution
$(14)(23)(56)$ of the $A_1$ fibers and a switch of the far leaves of
the $D_9$ fiber. Since the far leaves of the $D_9$ fiber are switched
by the Galois involution that multiplies the square root of $b'' =
I_{10}/4$ by $-1$, this tells us that we have the correct quadratic
twist, since $I_{10}$ is within a square factor of the discriminant of
the sextic. That is, making a different choice of $2$-level structure
while describing the fibration on the Kummer surface $Y$ would have
given us the same answer for $X$.

\section{An application to the moduli of six points in \texorpdfstring{$\Proj^1$}{P1}} \label{sec:6points}

The construction above gave us a correspondence of sextics
$$
f(x) = \sum f_i x^i  = f_6 \prod (x-x_i)
$$
and
$$
g(x) = \bigg(x^3 -\frac{I_4}{12}x + \frac{I_2 I_4 -3I_6}{108} \bigg)^2 + I_{10}\bigg(x - \frac{I_2}{24} \bigg).
$$
Therefore, over an algebraically closed field, we get a birational map
from the moduli space of $6$ points in $\Proj^1$ (i.e. the quotient of
$(\Proj^1)^6$ under the action of $\PGL2$ and $S_6$) with the space of
roots up to scaling of
$$
(x^3 + ax + b)^2 + (a'x+b')
$$
as $a,b,a',b'$ vary (we suppressed $b''$ since it just scales $a'$ and
$b'$).  This latter space is cut out inside $\Proj^5 =
\{(X_1:X_2:X_3:X_4:X_5:X_6)\}$ by the hyperplane $\sigma_1(X) = X_1 +
\ldots + X_6 = 0$ and the quartic hypersurface $\sigma_2(X)^2 = 4
\sigma_4(X)$, where $\sigma_2$ and $\sigma_4$ are the second and
fourth elementary symmetric functions of the $X_i$. Thus, we get a
model as a singular quartic threefold in $\Proj^4$, which is known in
the literature as the Igusa quartic. Here the Igusa quartic arises
naturally in the context of Shioda-Inose structures on K3 surfaces.

There is no simple one-one correspondence between the roots of $f(x)$
and $g(x)$, since the two actions of $S_6$ acting by the permutation
representation on the six roots of $f(x)$ on the six roots of $g(x)$
are related by an outer automorphism. To see this, we recall from the
last section that the permutation $(45)$ on the roots of $f(x)$ (the
Weierstrass points) acts on the roots of $g(x)$ (which are the
locations of the $A_1$ fibers) by the permutation $(14)(23)(56)$. By
symmetry, all the transpositions of $S_6$ act by a product of three
transpositions on the roots of $g(x)$. Thus we get a homomorphism
$S_6(f) \rightarrow S_6(g)$ which is an outer automorphism.

The curve $W$ defined by the equation $g(x) = 0$ is, as explained in
\cite{GLD}, the genus $2$ component of the fixed locus of the
involution $\tau$ which is associated to the K3 surface $X$ by virtue
of the fact that $\NS(X)$ has a $2$-elementary discriminant group
\cite{N2}. Explicitly, the involution $\tau$ is $(x,y,t) \mapsto
(x,-y,t)$ on the model of $X$ with $D_{14}$ fiber, an $A_1$ fiber and
a $2$-torsion section:
$$
y^2 = t^3 + \bigg(x^3 -\frac{I_4}{12} x + \frac{I_2 I_4 - 3I_6}{108} \bigg)t^2  -\frac{I_{10}}{4} \bigg(x- \frac{I_2}{24}  \bigg)t. 
$$

Consider the map $\phi$ from the moduli space ${\mathcal M}_2$ of
genus $2$ curves to itself which takes $C = V(f)$ to $W = V(g)$. That
is, we consider $g$ not only up to scaling of roots, but up to the
action of all elements of $\PGL2$ (on the roots) which fix the form $g
= p^2 + q$, where $p$ is a monic cubic polynomial with zero constant
term, and $q$ is a linear polynomial. Then, as noted in \cite{GLD},
$\phi$ is a rational self-map of ${\mathcal M}_2$ of degree $16$.

\section{Compatibility with formulas for \texorpdfstring{$\Km(E_1 \times E_2)$}{Km(E1xE2)} } \label{sec:product}

Here we check that the formulas we obtain degenerate to those of
Shioda \cite{S2} for the elliptic K3 surface corresponding to the
Kummer surface of a product of two elliptic curves, over an
algebraically closed base field. For the rest of this section, let us
fix this field $k$.

It is shown in \cite{S2} that given two elliptic curves $E_1$ and
$E_2$ with $j$-invariants $j_1$ and $j_2$, there is an elliptic K3
surface with two $E_8$ fibers which is $2$-isogenous to $\mathop{\rm
  Km}\nolimits(E_1 \times E_2)$. Its equation is given by
$$
y^2 = x^3 - 3\alpha x + \bigg(t + \frac{1}{t} - 2 \beta \bigg)
$$
with $\alpha = \sqrt[3]{j_1j_2}$ and $\beta =
\sqrt{(1-j_1)(1-j_2)}$. The choice of the square root and the cube
root here are arbitrary (since changing $x$ by a cube root of unity
multiplies $\alpha$ by it, whereas changing $y$ by a square root of
$-1$ and $x$ and $t$ by $-1$ changes the sign of $\beta$). Notice that
Shioda's $j$ is related to the $J$-invariant more commonly employed in
the literature by $j = J/1728$.

Now, the Igusa-Clebsch invariants $I_2, I_4, I_6, I_{10}$ (or
$A,B,C,D$ in Igusa's \cite{I2} notation) define a weighted projective
space. The $k$-points of this space with $I_{10} \neq 0$ correspond to
curves of genus $2$, therefore also their Jacobians with canonical
polarizations. Blowing up around the point $I_4 = I_6 = I_{10} = 0$
gives a variety which parametrizes Jacobians of genus $2$ curves as
well as their degenerations, products of two elliptic curves. Igusa
defines the coordinates
$$
x_1 = \frac{2^4 3^2 B}{A^2},\, x_2 = \frac{2^6 3^3 (3C-AB)}{A^3},\, x_3 = \frac{2^1 3^5 D}{A^5}
$$
$$
y_1 = \frac{x_1^3}{x_3},\, y_2 = \frac{x_2^2}{x_3},\, y_3 = \frac{x_1^2 x_2}{x_3}
$$

The point corresponding to a product of two elliptic curves with
$J$-invariants $J_1$ and $J_2$ has $x_1 = x_2 = x_3 = y_3 = 0$ and
$y_1 = J_1 J_2$, $y_2 = (J_1 - 1728)(J_2 - 1728)$.

Using these formulas Shioda's elliptic surface may be written as
$$
y^2 = x^3 - \frac{3}{12^2} \sqrt[3]{y_1} \, x + \bigg(t + \frac{1}{t} - \frac{2}{12^3} \sqrt{y_2}\bigg).
$$

Now let us transform the K3 surface
$$
y^2 = x^3 - t^3\bigg( \frac{I_4}{12} t +1\bigg)x + t^5\bigg(\frac{I_{10}}{4} t^2 + \frac{I_2 I_4 - 3 I_6}{108} t + \frac{I_2}{24}\bigg)
$$
by first setting $y = y t^6, x = x t^4$. This gives the equation
$$
y^2 = x^3 - \bigg( \frac{I_4}{12} + \frac{1}{t}\bigg)x + \bigg(\frac{I_{10} t}{4}  + \frac{I_2 I_4 - 3 I_6}{108} + \frac{I_2}{24 t}\bigg).
$$
Next, if we scale $x,y,t$ as $y = y/\lambda^3, x = x/\lambda^2, t = \mu t$, we get
$$
y^2 = x^3 - \lambda^4 \bigg( \frac{I_4}{12} + \frac{1}{\mu t}\bigg)x + \lambda^6 \bigg(\frac{I_{10} \mu t}{4}  + \frac{I_2 I_4 - 3 I_6}{108} + \frac{I_2}{24 \mu t}\bigg).
$$
Now we let
$$
\mu = \bigg(\frac{I_2}{6 I_{10}}\bigg)^{1/2} ,\, \lambda = \bigg(\frac{96}{I_2I_{10}}\bigg)^{1/12}.
$$
On simplification, the equation of the K3 surface becomes
$$
y^2 = x^3 - \bigg(\frac{3}{12^2} \sqrt[3]{\frac{2^{11} 3 B^3}{AD}} + \frac{1}{t}\bigg(\frac{2^2 D}{3^2 A^5}\bigg)^{1/6}  \bigg) x + \bigg(t + \frac{1}{t} + \frac{2}{12^3} \sqrt{\frac{2^{11} 3 (AB-C^2)}{AD}} \, \bigg)
$$
or
$$
y^2 = x^3 - \bigg(\frac{3}{12^2} \sqrt[3]{y_1} + \frac{1}{t}\bigg(\frac{2x_3}{3^7} \bigg)^{1/6} \bigg) x + \bigg(t + \frac{1}{t} + \frac{2}{12^3} \sqrt{y_2}\bigg).
$$
Now, noticing that $x_3 = 0$ for a product of two elliptic curves, we
get Shioda's formula, but with the sign of $\beta$ changed. However,
as we remarked earlier, this sign can be twisted away. Therefore, our
formulas indeed degenerate to those of Shioda.

\section*{Acknowledgements}

I thank Noam Elkies and Barry Mazur for suggesting the study of K3
surfaces as a thesis topic and for many insightful remarks and helpful
conversations. I also thank Tetsuji Shioda for discussions and
correspondence on elliptic surfaces, and Igor Dolgachev for helpful
remarks and pointing out some useful references. I also thank the
referees for their comments and corrections and Charles Doran for
bringing his work with Adrian Clingher to my attention.

\end{document}